\documentclass[10pt]{amsart}
\usepackage[utf8]{inputenc}
\usepackage{amssymb,latexsym}
\usepackage[mathscr]{eucal}
\usepackage{mathtools}

\theoremstyle{plain}
\newtheorem*{theorem}{Theorem}

\makeatletter
\@namedef{subjclassname@2020}{%
  \textup{2020} Mathematics Subject Classification}
\makeatother

\begin{document}

\title{An Unimaginative Proof of \\ Fermat's Two Squares Theorem}

\author[G. Bachman]{Gennady Bachman}
\address{Department of Mathematical Sciences\\ University of Nevada Las Vegas\\
4505 Maryland Parkway, Box 454020\\
Las Vegas, Nevada 89154-4020, USA}
\email{gennady.bachman@unlv.edu}

\date{March 2025}

\begin{abstract}
We give a simple direct proof of Fermat's two squares theorem. Our argument uses no intricate notions or ideas; one might say that it is a proof by careful bookkeeping. As such, the proof may be particularly easy to comprehend by students who are just starting out in math.
\end{abstract}



\maketitle

\begin{theorem} Every prime number $p$ congruent to 1 modulo 4 can be uniquely expressed as a sum of squares
\begin{equation}
p=a^2+b^2
\end{equation}
of two positive integers $a$ and $b$.
\end{theorem}

This is the Fermat's two squares theorem. It was conjectured by Girard in 1632, claimed by Fermat in 1659, with the first recorded proof of the existence of representation (1) given by Euler in 1749, and that of uniqueness of (1) by Gauss in 1801. The theorem is rich in history, old and new. The classic number theory textbook of Hardy and Wright \cite{HW} includes four different proofs. An extensive discussion of the history of this theorem and related topics will be found in Dickson \cite[Volume 2, Chapter VI]{Di}. For a reader interested in just the historical highlights of this theorem we highly recommend the introduction to the paper \cite{CELV} of Clarke et al. The current developments connected to this result were sparked by the papers of Heath-Brown \cite{HB} and Zagier \cite{Za}. The recent papers of Elsholtz \cite{El} and Bacher \cite{Ba} contain some of the newest developments as well as an engaging discussion of the history of this theorem up to the present.

It is rather inconceivable that the argument we present in this note is new and somehow evaded detection all these years. But the author failed to find it in the literature and this note is the toll for that failure. In addition to being particularly simple and direct, our argument establishes existence and uniqueness in one go (as opposed to most other proofs). Our starting point is that $-1$ is a quadratic residue for every prime $p\equiv1\pmod 4$ and that (1) would not have been possible otherwise. But in keeping with our slogan ``simple and direct'' and to keep this note self-contained, we include a quick verification of this fact. The set of integers $\{2,3,\dots,p-2\}$ partitions into an odd number of pairs $\{u,v\}$ with $uv\equiv1\pmod p$. It follows that it cannot be partitioned into 2-element parts $\{u',v'\}\,(=\{p-u,v\})$ with $u'v'\equiv-1\pmod p$. In other words, there must be a pair $\{u,v\}$ with $v=p-u$, or $u^2\equiv-1\pmod p$.

For a given prime $p\equiv1\pmod4$, we fix $u$ such that $u^2\equiv-1\pmod p$ and put $R=\lfloor\sqrt p\rfloor$. Then (1) implies that either $a\equiv ub\pmod p$ or $b\equiv ua\pmod p$, and to find solutions to this equation is to find $1\le a,b\le R$ with $b\equiv ua\pmod p$, say. It is then natural to associate with each integer $x$ in the range $1\le x\le R$, the corresponding $x\mapsto y(x)$ with $1\le y(x)<p$ and $y(x)\equiv ux\pmod p$. This gives $R$ distinct values $y(x)$ and we arrange them in increasing order
\[ 0<y_1<y_2<\dots<y_R<p. \]
This arrangement induces the associated arrangement $\bigl(x_i\bigr)_{i=1}^R$ of integers $1\le x\le R$, where $ux_i\equiv y_i\pmod p$. The pairs $(x_i,y_i)$ satisfy $x_i^2+y_i^2\equiv0\pmod p$, and the equation $x_i^2+y_i^2=p$ holds for every $y_i\le R$. It is now immediate that the claim of the theorem is equivalent to the assertion that $y_1\le R<y_2$.

First, we consider the possibility that $y_1<y_2\le R$, so that
\[ x_1^2+y_1^2=x_2^2+y_2^2=p \quad\text{and}\quad x_2<x_1. \]
Then $0<(y_2-y_1), (x_1-x_2)<R$ and $x_1-x_2\equiv u(y_2-y_1)\pmod p$. It follows that $(y_2-y_1)^2+(x_1-x_2)^2=p$ is another solution. But this is absurd since the left side $=2p-2(x_1x_2+y_1y_2)$ is even.

It remains to show that $y_1\le R$. Observe that 
\[(y_1-0)+(y_2-y_1)+\dots+(y_R-y_{R-1})+(p-y_R)=p,\]
so one of these $R+1$ differences must be $<R$. If $p-y_R<R$, then we interpret $p-y_R$ as one of the $x_j$s and then $x_R$ as the corresponding $y_j$:
\[ x_R\equiv u(p-y_R)\equiv ux_j\equiv y_j\pmod p. \]
This says that $y_1\le x_R\le R$. Next is the possibility that $y_{i+1}-y_i<R$ for some $i$. Now, if $x_{i+1}>x_i$, then $x_{i+1}-x_i=x_j$ and $y_{i+1}-y_i=y_j$ for some $j$:
\[ y_{i+1}-y_i\equiv u(x_{i+1}-x_i)\equiv ux_j\equiv y_j\pmod p, \]
yielding $y_1\le y_{i+1}-y_i<R$. If, on the other hand, $x_{i+1}<x_i$, then we interpret $y_{i+1}-y_i$ as one of the $x_j$s and then $x_i-x_{i+1}$ as the corresponding $y_j$:
\[ x_i-x_{i+1}\equiv u(y_{i+1}-y_i)\equiv ux_j\equiv y_j\pmod p, \]
concluding that $y_1\le x_i-x_{i+1}<R$ again. We have thus verified that $y_1\le R$, as required.

We conclude this note with two remarks. In his ``Apology'' \cite{Ha} Hardy commented: ``This is Fermat's theorem, which is ranked, very justly, as one of the finest in arithmetic. Unfortunately there is no proof within the comprehension of anybody but a fairly expert mathematician.'' One dares to imagine that Hardy would have approved of our argument as having its place amongst other proofs of this theorem.

For our second remark, we introduce the notation $\langle n\rangle_m$ to denote the least nonnegative residue of $n$ modulo $m$:
\[\langle n\rangle_m\equiv n\pmod m \quad\text{and}\quad 0\le\langle n\rangle_m<m. \] 
With this notation our function $y(x)$ takes the form $y(x)=\langle ux\rangle_p$, and the main part of our argument showed that
\[ \min_{1\le x<\sqrt p}\langle ux\rangle_p<\sqrt p. \]
Exploring this notion further leads to a slightly weaker inequality but a more general result
\[ \min_{1\le|x|\le\sqrt m}\langle ax\rangle_m<\sqrt m, \]
valid for all $m\ge2$ and $(a,m)=1$. This result is due to Thue \cite{Th} (restated here in the form matching the forgoing discussion). This more insightful approach leads to a more elegant proof of Fermat's theorem \cite[Lemma 2.13]{NZM}. Still, it is fun to puzzle over the fact that the argument of the form presented here was not one of the earlier proofs of the theorem.

\end{document}